\documentstyle[12pt]{article}
\def\Bbb{{}}

\newcommand{\la}{\langle}
\newcommand{\ra}{\rangle}

\newcommand{\mod}{\mbox{mod\,\,}}

\newcommand{\Ind}{\mbox{Ind}}

\newcommand{\f}{\,\,\,\forall\,\,\,}

\newcommand{\supp}{\mbox{supp}}

\newcommand{\LL}{{\cal L}}

\def \N{\hbox{$I\hskip -4pt N$}}
\def \Z {\hbox{$Z\hskip -5.2pt Z$}}
\def\sZ{\hbox{$\sc Z\hskip -4.2pt Z$}}

\def \C{\hbox{$C \hskip -7pt \vrule height 6pt depth 0pt \hskip 6pt$}}

\def\qed{\hfill \hfill \ifhmode\unskip\nobreak\fi\ifmmode\ifinner
         \else\hskip5pt\fi\fi
 \hbox{\hskip5pt\vrule width4pt height6pt depth1.5pt\hskip 1 pt}}
\def\a{\alpha}
\def\b{\beta}
\def\d{\delta}

\def\l{\lambda}

\def\o{\omiga}

\def\Vir{{\rm{Vir}}}

\def\sc{\scriptstyle}

\def\cl{\centerline}

\def \LL{{\cal L}}

\def \o{\otimes}

\def\sz{{Z\hskip -3.8pt Z}}
\def \B{{\cal B}}
\def \ht{\hbox{ht}}
\def \SS{\cal S}
\def \hm{\hbox{hm}}
\begin{document}
\par\
\par\
\par

\cl{\large{\bf  Classification of irreducible Harish-Chandra
modules }}

\cl{\large{\bf over the loop-Virasoro algebra}\footnote{AMS
Subject Classification: 17B10,  17B65, 17B68.\\
\indent \hskip .3cm  Research supported by  NSERC, Postdoctoral
Research Grant, and NSF of China.\\
\indent \hskip .3cm Keywords: loop-Virasoro algebra, weight
module, Harish-Chandra module}}

\par
\vskip 15pt
 \centerline{Xiangqian
Guo, Rencai Lu and Kaiming Zhao}

\par\vskip 10pt
\begin{abstract}
{The loop-Virasoro algebra is the Lie algebra of the tensor
product of the Virasoro algebra and the Laurent polynomial
algebra. This paper classifies  irreducible Harish-Chandra modules
over the loop-Virasoro algebra, which turn out to be highest
weight modules, lowest weight modules and evaluation modules of
the intermediate series (all wight spaces are 1-dimensional). As a
by-product, we obtain a classification of irreducible
Harish-Chandra modules over truncated Virasoro algebras.

We also determine the necessary and sufficient conditions for
highest weigh irreducible modules over the loop-Virasoro algebra
to have all finite dimensional weight spaces, as well as the
necessary and sufficient conditions for highest weigh Verma
modules  to be irreducible.}
\end{abstract}

\vskip .5cm

\par
\cl{{\bf \S1. Introduction}}
\par
\vskip .3cm Representations of many infinite dimensional Lie
algebras have important applications in mathematics and physics.
The relation to physics is well established in the book on
conformal field theory [FMS]. Recently, developing representation
theory for various infinite dimensional Lie algebras has attracted
extensive attention of many mathematicians and physicists. Many of
these Lie algebras are related to the Virasoro algebra.

\vskip 5pt The Virasoro algebra theory has been widely used in
many physics areas and other mathematical branches, for example,
 conformal field theory [IKUX, IUK], combinatorics  [M], Kac-Moody algebras [K, MoP], vertex
algebras [DL], and so on. Let us first recall the definition of
this Lie algebra. \vskip 5pt

In this paper we denote by $\C$,  $\Z$, $\Z_+$ and $ \N$ the set
of complex numbers, integers, nonnegative integers and positive
integers respectively.

\vskip 5pt The {\bf Virasoro algebra $\Vir:=\Vir[\Z]$} (over $\C$)
is  the Lie algebra with the basis
$\bigl\{c,d_{i}\bigm|i\in\Z\bigr\}$ and the Lie brackets defined
by

$[d_m,d_n]=(n-m)d_{m+n}+\delta_{m,-n}\frac{m^{3}-m}{12}c,\qquad\forall
m,n\in \Z,$

$[d_m,c]=0, \qquad\forall m\in \Z.$ \vskip 5pt

The structure theory of Harish-Chandra modules over the Virasoro
algebra has been  developed fairly well. For details, we refer
readers to [Ma], [MP], [MZ], the book [KR] and  references
therein. In particular, the Kac conjecture, i.e., the
classification of irreducible Harish-Chandra modules over the
Virasoro algebra was obtained in [Ma].

\vskip 5pt It is also well-known that simply-laced affine
Kac-Moody algebras are loop algebras of finite dimensional simple
Lie algebras, and their representation theory is very different
from that of finite-dimensional Lie algebras, much richer and
having more applications to other fields [K]. This stimulates us
to investigate the loop algebra of the Virasoro algebra, which we
simply call the loop-Virasoro algebra. Other reasons for us to
study the loop-Virasoro algebra are the paper [Wi] on truncated
Virasoro algebras $\Vir\otimes \Big(\C[t]/t^n\C[t]\Big)$ and the
paper [ZD] on a special truncated Virasoro algebra $W(2,2)$. We
hope this theory will also have applications to physics and other
areas in mathematics. Let us define the loop-Virasor algebra
precisely.

 \vskip 5pt
The {\bf loop-Virasoro algebra} $\LL$ is the Lie algebra that is
the tensor product of the Virasoro Lie algebra $\Vir$ and the
Laurent polynomial algebra $\C[t^{\pm1}]$, i.e., $\LL=\Vir\otimes
\C[t^{\pm1}]$ with a basis $\{ c\otimes t^j, d_i\otimes t^j |
i,j\in \Z \}$ subject to the commutator relations:
$$[d_i\otimes t^j,d_k\otimes t^l]=(k-i)\left(d_{i+k}\otimes
t^{j+l}\right)+\delta_{i+k,0}\frac{i^3-i}{12}\left(c\otimes
t^{j+l}\right),$$
$$[d_i\otimes t^j,c\otimes t^l]=0.$$

We see that $\LL$ has a copy of $\Vir$ which is $\Vir\otimes 1$.
For simplicity, we shall write $d_i=d_i\otimes 1$, $c=c\otimes 1$,
$d_i(j)=d_i\otimes t^j$ and $c(j)=c\otimes t^j$.

Let  $\LL_i=(\C d_i\oplus \d_{i,0}\C c)\otimes \C[t^{\pm}]$. Then
$\LL$ has the natural $\Z$-gradation $\LL=\oplus_{i\in\sz}\LL_i$.
If we denote $\LL_+=\bigoplus_{i>0}\LL_i$ and
$\LL_-=\bigoplus_{i<0}\LL_i$, then $\LL$ has the triangular
decomposition $$\LL=\LL_-\oplus\LL_0 \oplus\LL_+.\eqno(1.1)$$

Note that $\LL_0$ is an infinite dimensional abelian subalgebra of
$\LL$, and that $c\otimes\C[t^{\pm1}]$ is the center of $\LL$. The
loop-Virasoro algebra $\LL$ is not a pre-exp-polynomial algebra
defined in Definition 1.5 [BGLZ] since Condition (P1) cannot be
satisfied. It still has similar properties to pre-exp-polynomial
algebras.

\vskip 5pt The paper is organized as follows. In Section 2,  we
recall some related concepts, results for later use, and define
the highest weight Verma $\LL$-module $\bar V(\varphi)$ and the
corresponding highest weight irreducible $\LL$-module $V(\varphi)$
for any linear map $\varphi:\,\LL_0\to\C$. In Section 3, using
Matieu's Theorem we give a rough classification of irreducible
Harish-Chandra modules over $\LL$, which turn out to be highest
weight modules, lowest weight modules, and uniformly bounded
modules (Theorem 3.1). In Section 4, using Jacobson radical theory
on finite dimensional associative algebras we prove that uniformly
bounded modules over the loop-Virasoro algebra are of the
intermediate series, i.e., all weight spaces are $1$-dimensional.
Consequently, irreducible Harish-Chandra modules over the
loop-Virasoro algebra are highest weight modules, lowest weight
modules and  modules of the intermediate series. As a by-product,
we also obtain a classification of irreducible Harish-Chandra
modules over truncated Virasoro algebras studied in [Wi], as well
as a classification of irreducible Harish-Chandra modules over the
algebra $W(2,2)$ studied in [ZD]. In Section 5, we describe weight
modules over $\LL$ with all weight spaces $1$-dimensional, which
we prove to be evaluation modules of the intermediate series
(Theorem 5.1). In Section 6, we obtain the necessary and
sufficient conditions for $V(\varphi)$ to have all finite
dimensional weight spaces (Theorems 6.1, 6.6), as well as the
necessary and sufficient conditions for $\bar V(\varphi)$ to be
irreducible (Theorems 6.5, 6.6). Theorem 6.4 reduces the study of
highest weight irreducible modules $V(\varphi)$ with all finite
dimensional weight spaces to the study of such modules over
truncated Virasoro algebras $\Vir\otimes \Big(\C[t]/t^n\C[t]\Big)$
which were investigated in [Wi].  \vskip 15pt

\par
\cl{{\bf \S2. Modules over $\Vir$}}
\par
\vskip .2cm In this section, by a module it means a module over
$\Vir$, $\LL$ or $\Vir\otimes \Big(\C[t]/t^n\C[t]\Big)$.

A module $V$ is called {\bf trivial} the action of the Lie algebra
is zero. For any module $V$ and $\lambda\in \C$, the subspace
$V_{\lambda}:=\bigl\{v\in V\bigm|d_{0}v=\lambda v\bigr\}$ is
called the {\bf weight space} of $V$ corresponding to the weight
$\lambda$. A module $V$ is called a {\bf weight module} if $V$ is
the sum of its weight spaces, and a weight module is called a {\bf
Harish-Chandra module} if all its weight spaces are finite
dimensional.

For a weight module $V$, we define $\supp V:=\bigl\{\lambda\in
\C\bigm|V_{\lambda}\neq 0\bigr\}$, which is generally called the
{\bf weight set} (or the {\bf support}) of $V$. A weight module
$V$ is said to be {\bf uniformly bounded}, if there exists $N\in
\N$ such that $\dim V_{x}<N$ for all $x\in \supp V$.

Let $U:=U(\Vir)$ be the universal enveloping algebra of the
Virasoro algebra $\Vir$. For any $\dot c, h\in F$, let $I(\dot
c,h)$ be the left ideal of $U$ generated by the set $$
\bigl\{d_{i}\bigm|i>0\bigr\}\bigcup\bigl\{d_0-h\cdot   1, c-\dot
c\cdot 1\bigr\}. $$ The {\bf Verma module} with highest weight
$(\dot c, h)$ for $\Vir$ is defined as $\bar M(\dot c,h):=U/I(\dot
c,h)$. It is a highest weight module of $\Vir$ and has a basis
consisting of all vectors of the form $$ d_{-i_1}d_{-i_2}\cdots
d_{-i_k}v_{h};\quad k\in{\N}\cup\{0\}, i_{j}\in\N,
i_{k}\geq\cdots\geq i_2\geq i_1>0.
$$ Then we have the irreducible quotient $ M(\dot c,h)=\bar M(\dot c,h)/J$
where $J$ is the maximal proper submodule of $\bar M(\dot c,h)$.
It is well known that, if $(\dot c,h)\ne (0,0)$, for any $N\in\N$
there exists $k\in\N$ such that $\dim M(\dot c,h)_{h-k}>N$. (See,
for example, Claim 4 on Page 647 in [LZ2]).

It is well known that [SZ] a $\Vir$-module of the intermediate
series must be one of $V(\alpha,\beta),A_{a},B_{b}$ for some
$\alpha,\beta,a,b\in \C$ or one of their quotient submodules,
where $V(\alpha,\beta)$ (resp. $A_{a}$, $B_{b}$) has basis
$\{v_{\alpha+i}|i\in \Z\}$ (resp. $\{v_{i}|i\in \Z\}$) with
trivial central actions and

\smallskip
$V(\alpha,\beta):\ \ \
d_{i}v_{\alpha+k}=(\alpha+k+i\beta)v_{\alpha+k+i}$,
\smallskip

$A_{a}:\ \ \ d_{i}v_{k}=(k+i)v_{k+i}, \forall k\neq 0;\ \
d_{i}v_{0}=i(i+a)v_{i}$,
\smallskip

$B_{b}:\ \ \ d_{i}v_{k}=kv_{k+i}, \forall k+i\neq 0;\ \
d_{i}v_{-i}=-i(i+b)v_{0}$.
\smallskip

The bases given above are called {\bf the standard bases} of the
corresponding modules. It is known that $V(\alpha,\beta)\cong
V(\alpha+1,\beta), \forall \alpha,\beta$ and that $V(\alpha,0)\cong
V(\alpha,1)$ if $\alpha\notin\Z$. It is also clear that $V(0,0)$ and
$B_{b}$ both have $\C v_0$ as a submodule, and their corresponding
quotients are isomorphic, which we denote by $V'(0,0)$. Dually,
$V(0,1)$ and $A_{a}$ both have $\C v_0$ as a quotient module, and
their corresponding submodules are isomorphic to $V'(0,0)$. For
convenience we simply write $V'(a,b)=V(a,b)$ when $V(a,b)$ is
irreducible.

Let us consider the graded dual modules of the above modules.
Using the standard basis for $V=A_a$, we define $v^*_i\in V^*$
such that $v^*_i(v_{-j})=\delta_{i,j}$. As a $\Vir$-module, it is
easy to verify that the dual module $V^*=\oplus_{i\in\Z}\C v^*_i$
is isomorphic to $B_a$ with standard basis $\{v^*_i|i\in\Z\}$,
i.e., $A_a^*\simeq B_a$. Similarly, $V(a, b)^*\simeq V(a, 1-b)$.
\medskip

 {\bf Modules of the intermediate series}
 are  indecomposable weight modules with all weight
spaces $1$-dimensional.

For any $\Vir$-module of intermediate series $V$, we can define the
{\bf evaluation module} $V(e)$ for any nonzero $e\in\C$ over $\LL$
as follows: $V(e)=V$ as vector spaces and the actions of $\LL$ on
$V(e)$ is given by $(d_i\o t^m)u=e^md_iu, \f u\in V$.

 The evaluation module $ M(\dot c,h)(
a)$ over  $\LL$   were studied in [W]. We will see in Section 7
that these irreducible modules are only a small portion of highest
weight irreducible modules over $\LL$.

\medskip
Similar to the definition in [RZ] or in [BGLZ], now we define
highest weight modules over $\LL$. For any linear map
$\varphi:\LL_0\to\C$ with $\varphi(d_0)=\lambda$, we define the
1-dimensional $(\LL_0+ \LL_+)$-module $\C v_0$ via
 $$\LL_iv_0=0,\,\, {\rm if} \,\,i>0; \,\,xv_0=\varphi(x)v_0,\f x\in L_0.
 \eqno(2.1)$$
 Consider the induced $\LL$-module
 $$\bar V(\varphi)=\Ind^{\LL}_{\LL_++\LL_0}\C v_0=U(\LL)\otimes_{U(\LL_++\LL_0)}\C v_0,$$
 where $U(\LL)$ is the universal enveloping algebra of  $\LL$.
 It is clear that, as vector spaces,
 $\bar V(\varphi)\simeq U(\LL_-)$. The  module $\bar V(\varphi)$
 has a unique maximal proper submodule $J$. Then we obtain the irreducible
 module
 $$V(\varphi)={\bar V(\varphi)\over J}.\eqno(2.2)$$
 This module is called a {\bf highest weight module over} $\LL$.
 We see that $V(\varphi)$ is uniquely determined by the linear
 function $\varphi$, and that $V(\varphi)=\oplus_{i\in\sz_+}V_{\lambda-i}$ where
 $$V_{\lambda-i}=\{v\in V(\varphi)\,\,|\,\,d_0v=(\lambda-i)v\}.\eqno(2.3)$$

 Similarly we can define highest weigh modules over truncated
 Virasoro algebras $\Vir\otimes\Big(\C[t]/\langle t^n\rangle\Big)$ for
 any $n\in\N$.

 \vskip .5cm
\par
\cl{{\bf \S3. Classification of Harish-Chandra modules over
$\LL$}}
\par

In this section, we shall give a rough classification of irreducible
Harish-Chandra modules over  $\LL$, i.e., we shall prove the
following

{\bf Theorem 3.1.} {\it  Let $V$ be an irreducible Harish-Chandra
$\LL-$module. Then $V$ is either a highest weight module,  a
lowest weight module or a uniformly bounded module.}

{\it Proof.} Assume that $V$ is not uniformly bounded. As a $\Vir$
module, suppose that $W$ is a minimal $\Vir$ submodule of $V$ such
that $V/W$ is trivial and that $T$ is the largest trivial $\Vir$
submodule of $W$.

The result in [MP] ensures that the $\Vir$-submodule decomposition
$\bar{W}:=W/(W\cap T)\cong \bar{W}^+\oplus \bar{W}^-\oplus
\bar{W}^0$ where $\supp(\bar{W}^+)$   is upper bounded,
$\supp(\bar{W}^-)$ is lower bounded and $\bar{W}^0$ is uniformly
bounded. Without loss of generality, we may assume that
$\bar{W}^+$ is nontrivial. For any $w\in W$, we denote by
$\bar{w}$ the image of $w$ in $\bar{W}$.

Since the action of $c$ on $V$ is a scalar, we consider two cases
according to the action of $c$.

{\bf Case 1.} The action of $c$ is zero.

Let $\mu$ be the highest weight of $\bar{W}^+$. Take any $v_1\in
{W}^+_{\mu}$ such that $\bar v_1\ne0$.

In case that $\mu\neq 0$, we set $\l=\mu$ and that $v_0=v_1$.

In case that $\mu=0 $, we consider the $\Vir$ module $U(\Vir)v_1$.
Because of the choice of $W$,  $U(\Vir)v_1/(U(\Vir)v_1\cap T)$ is a
highest weight $\Vir$ submodule of $\bar{W}^+$ and
$U(\Vir)v_1=U(\Vir_-)v_1+(U(\Vir)v_1\cap T)$ as vector spaces. It is
easy to see that there exists a weight vector $v_0\in U(\Vir)v_1$
with  weight $-1$ (i.e., $d_0v_0=-v_0$) such that $\bar{v_0}$ is a
highest weight vector in $U(\Vir)v_1/(U(\Vir)v_1\cap T)$, i.e.,
$d_iv_0=0$ for all $i\in\N$. Let $\l=-1$ be the weight of $v_0$.

Now in both cases $\l\neq 0$. Let $M=U(\Vir)v_0$ and $T'=T\cap M$.
Similarly we have that $M/T'\subset \bar{W}^+$ is a nontrivial
highest weight $\Vir$ module and that $M=U(\Vir_-)v_0+T'$ as
vector spaces. Let $M'$ be the largest $\Vir$ submodule of $M$
with $M'_{\l}=0$. Then $T'\subset M'$, and $M/M'$ is isomorphic to
the nontrivial irreducible $\Vir$ module $M(0,\l)$.

There is some $k\in \N$ such that $\dim M(0,\l)_{\l-k}>\dim V_{\l}
$. For any $i\in \Z$, consider the linear map:
$$d_k\o t^i: M''_{\l-k} \to V_{\l},$$
where $M''_{\l-k}$ is a complement of $M'_{\l-k}$ in $M_{\l-k}$.
Since $\dim M''_{\l-k}\ge \dim M(0,\l)_{\l-k}>\dim V_{\l} $  there
is some $w_i\in M''_{\l-k}\setminus\{0\}\subset M_{\l-k}\setminus
M'_{\l-k}$ such that $(d_k\o t^i)w_i=0$.

Choose $k_0\in \N$ such that $\l+j\neq 0$ for any $j\in \N$ with
$j>k_0$. For any $j>k_0$, noticing that $(M'/T')_{\l+j}=0$,
$(M/M')_{\l+j}=0$, and $d_{k+j}w_i\in M_{\l+j}$, then $M_{\l+j}=0$
and $d_{k+j}w_i=0$ for any $j>k_0$ and $i\in\Z$. Thus we have that
$$0=d_{k+j}(d_k\o t^i)w_i=-j(d_{2k+j}\o t^i) w_i,\forall j>k_0.$$

Since $w_i\in M_{\l-k}\setminus M'_{\l-k}$, then there are some
$d_{i_1},...,d_{i_r}\in \Vir$ with $i_1,...i_r>0, i_1+...+i_r=k$,
such that $0\neq d_{i_1}...d_{i_r}\bar{w_i}\in
(M/M')_{\l}=M_{\l}/M'_{\l}$. since $M'_{\l}=0$ and $\l\neq 0$, we
have that $0\neq d_{i_1}...d_{i_r}w_i\in M_{\l}=\C v_0$. Thus
there is some $ a\in \C$ such that $ad_{i_1}...d_{i_r}w_i=v_0$.

Using the fact that  $(d_{2k+j}\o t^i)w_i=0$ for all $j>k_0$, we
deduce that $(d_{2k+j}\o t^i)v_0=(d_{2k+j}\o
t^i)(ad_{i_1}...d_{i_r})w_i=0$, for any $j>k_0, i\in \Z$.

Let $V^+=\{v\in V | \mbox{ there is some } n\in \N \mbox{ such
that } (d_{j}\o t^i)v=0, \mbox{ for any } j>n, i\in \Z \}.$
Clearly,  $V^+\neq 0$, and it is straightforward to check that
$V^+$ is an $\LL$-submodule of $V$. Thus $V=V^+$.

 If $\bar{W}^-\ne0 $ or $\bar{W}^0\ne0 $, by Lemma 1.6 in [Ma] or the structure
 of $\bar{W}^ 0 $ there
exists an element $v$ in $W$ which cannot be annihilated by $d_k$
for all sufficiently large $k$. So $v\notin   V^+$, a
contradiction. We obtain that $\bar{W}^-=0$ and $\bar{W}^0 =0$.
Then $\supp(V)$ itself is upper bounded. Therefore $V$ is a
highest weight $\LL$-module.

{\bf Case 2.} The action of $c$ is nonzero.

This case is much simpler. Since $c$ acts as a nonzero scalar,
 $V$ does not possess any trivial sub-quotient $\Vir$ modules.
Then in the argument of Case 1, $V=W=\bar{W}$ and $\bar{W}^0=0$.
The argument of the proof for Case 1 with minor modifications goes
more smoothly. \qed

Applying Theorem 3.1 to the the truncated Virasoro algebra
$\Vir\otimes \Big(\C[t]/t^n\C[t]\Big)\simeq \Vir\otimes
\Big(\C[t^{\pm1}]/(t+1)^n\C[t^{\pm1}]\Big)$ we can classify
irreducible Harish-Chandra module over the truncated Virasoro
algebra. (See [Wi] for the definitions).

 {\bf Corollary 3.2.} {\it Let $V$ be an irreducible Harish-Chandra
module over the truncated Virasoro algebra $\Vir\otimes
\C[t]/\langle t^n\rangle$ for any $n\in\N$. Then $V$ is either a
highest weight module, a lowest weight module or a uniformly
bounded module.}

 \vskip .5cm
\par
\cl{{\bf \S4. Complete classification of Harish-Chandra modules
over $\LL$}}
\par

In this section we shall give a complete classification of
irreducible Harish-Chandra  modules over $\LL$. We first prove the
result for the truncated Virasoro algebra $G=\Vir\otimes
\Big(\C[t]/\langle t^2\rangle\Big)$. In $G$, we shall write
$d_i=d_i\otimes 1$, $c=c\otimes 1$, $d'_i=d_i\otimes t$ and
$c'=c\otimes t$. Let $H$ be the abelian subalgebra of $G$ spanned by
$d'_i,c'$ $(i\in\Z)$.

{\bf Theorem 4.1.} {\it  Let $V$ be a uniformly bounded irreducible
$G-$module. Then $HV=0$.  Consequently, $V$ is an irreducible weight
module over $\Vir=\oplus_{i\in\sZ}\C d_i\oplus\C c$.}

{\it Proof.} We first introduce some notations. Denote by $U(G)$
and $U(H)$ the universal enveloping algebras of $G$ and $H$
respectively. Note that $U(H)$ is a commutative associative
polynomial algebra in the variables: $d'_i,c$ ($i\in\Z$). Let
$U(G)_i=\{x\in U(G)\,|\,[d_0, x]=ix\}$ and $U(H)_i=\{x\in
U(H)\,|\,[d_0, x]=ix\}$ for any $i\in\Z$. Then $U(G)_0$ and
$U(H)_0$ are associative algebras and $U(H)_0\subset U(G)_0$. Then
${\cal H}=U(G)_0\cap \Big( U(G)HU(G)\Big)$ is an ideal of
$U(G)_0$.

 Consider the weight space decomposition
$V=\oplus_{i\in\sZ}V_{\a+i}$ where $\a\in\C$. We may assume that
$\dim V_{\a+i}\le n$ for some $n\in \N$. We know that
$\LL_iV_{\a+j}\subset V_{\a+i+j}$ for any $i,j\in\Z$. Then there
exists a nonzero homogeneous vector $v\in V_{\a+i_0}$ and $a\in\C$
such that $d'_0v=av$, i.e., $(d'_0-a)v=0$. For any fixed $i\in\Z$
there exist monomials $x_1, x_2,\cdots,x_n\in U(G)_{i-i_0}$ of the
form $x_k=d_{j_1}d_{j_2}\cdots d_{j_s}x'_k$ where $j_1\le j_2\le
\cdots\le j_s$ ($s$ depends on $k$) and $x'_k\in U(H)$, such that
$V_{\a+i}={\rm span}\{x_kv:1\le k\le n\}$. It is easy to check that
$(d'_0-a)^{s+1}x_kv=0$. Thus $(d'_0-a)$ is nilpotent on each
$V_{\a+i}$. Since $\dim V_{\a+i}\le n$ we know that $(d'_0-a)^n=0$
on each $V_{\a+i}$. Consequently $(d'_0-a)^nV=0$.

Computing $d_j^r(d'_0-a)^nV=0$ we deduce that
$(d'_j)^r(d'_0-a)^{n-r}V=0$ for any $j\ne 0$.

{\bf Claim 1.} $a=0$.

To the contrary, suppose $a\ne0$. Since $c'$ acts as a scalar on
$V$. we still denote this scalar as $c'$. Choose $i_0\in\N$ such
that $2i_0a+\frac{i_0^3-i_0}{12}c'\ne0$. Suppose $m$ is minimal
such that $(d'_{i_0})^{m}V=0$. Since $V\ne0$,  $m>0$. Computing
$d_{-i_0}(d'_{i_0})^mV=0$ we deduce that
$(d'_{i_0})^{m-1}(2i_0d'_0+\frac{i_0^3-i_0}{12}c')V=0$. Combining
with $(d'_{i_0})^{m-1}(d'_0-a)^{n-m+1}V=0$ we obtain
$(d'_{i_0})^{m-1}V=0$ which contradicts the choice of $m$. Thus
Claim 1 follows.

From the proof of Claim 1 we also see that $c'V=0$.

{\bf Claim 2.} For any $i_1,i_2,\cdots,i_n\in\Z$,
$d'_{i_1}d'_{i_2}\cdots d'_{i_n}V=0$.

We first show by induction on $r: 0\le r\le n$ that
$d'_{i_1}d'_{i_2}\cdots d'_{i_r}(d'_0)^{n-r}V=0$ For any
$i_1,i_2,\cdots,i_r\in\Z$. The statement for $r=0$ is trivial.
Suppose it is true for $r (<n)$. Applying $d_{i_{r+1}}$ to
$d'_{i_1}d'_{i_2}\cdots d'_{i_r}(d'_0)^{n-r}V=0$ and using the
inductive hypothesis we deduce that $d'_{i_1}d'_{i_2}\cdots
d'_{i_{r+1}}(d'_0)^{n-r-1}V=0$. Claim 2 follows.

From Claim 2 we see that $H^nV=0$. Then $U(G)H^nU(G)V=0$. It is
not hard to verify that $(U(G)HU(G))^n=U(G)H^nU(G)$. Then ${\cal
H}^nV=0$.

Let $f$ be the representation of $U(G)_0$ on $V_{\a+i}$ for an
arbitrary  fixed $i\in\Z$. Since $V_{\a+i}$ is a
finite-dimensional simple faithful module over the associative
algebra $U(G)_0/\ker(f)$,  $U(G)_0/\ker(f)$ is a
finite-dimensional simple associative algebra and its Jacobson
radical $J=0$. Since $\Big(({\cal H}+\ker(f))/\ker(f)\Big)^n$
$V_{\a+i}=0$, we have $\Big(({\cal H}+\ker(f))/\ker(f)\Big)^n=0$.
Hence $({\cal H}+\ker(f))/\ker(f)\subset J=0$ (see [Lemma 9.7.2,
A]). Therefore $({\cal H}+\ker(f))/\ker(f)=0$, i.e., ${\cal H}\in
\ker(f)$. In particular, $d'_0\in{\cal H}$, that is, $d_0V=0$.

For any $k\in\Z$, since $kd'_kV=[d_0'd_k]V=0$, we know
$HV=0$.\qed

{\bf Theorem 4.2.} {\it Let $\LL(n)= \Vir\otimes (\C[t]/\langle
t^n\rangle)$ where $ n\geq 2$. Suppose that $V$ is a nontrivial
uniformly bounded irreducible $\LL(n)$ module. Then $(\Vir\otimes
t\C[t])V=0$ and $V$ is an irreducible weight module over
$\Vir=\oplus_{i\in\sZ}\C d_i\oplus\C c$.}

{\it Proof}. We will show this theorem by induction on $n$. The
result for $n=2$  is Theorem 4.1. Now we suppose that the theorem is
true for $n$ $(\ge 2$) and that $V$ is a uniformly bounded
irreducible $\LL(n+1)$ module. Let $L'=\Vir\oplus (\Vir\otimes
t^n)$. Then $V$ can be viewed as an $L'$ module, which is isomorphic
to $ \LL(2)\simeq G$. From Theorem 4.1 we know that $V$ has an
irreducible $L'$ submodule $V'$ on which $\Vir\otimes t^n$ act
trivially. Hence we have that $V=U(\Vir\otimes 1)U(\Vir\otimes
t^1)\cdot U(\Vir\otimes t^{n-1})V'$. It is easy to verify that
$\Vir\otimes t^n$ acts trivially on $V$. Thus $V$ is a module over
$\LL(n)\simeq \LL(n+1)/( \Vir\otimes t^n)$ (as a natural quotient
algebra of $\LL(n+1)$). By the induction hypothesis, we know that
$(\Vir\otimes t\C[t])V=0$ and $V$ is an irreducible weight module
over $\Vir=\oplus_{i\in\sZ}\C d_i\oplus\C c$.\qed

{\bf Theorem 4.3.} {\it Let $f(t)\in\C[t]$ have a nonzero constant
term and be of degree $\ge1$, $\Bbb L(f)=\Vir\otimes
(\C[t^{\pm1}]/\langle f\rangle) $. If $V$ is a nontrivial uniformly
bounded irreducible $\Bbb L(f)$ module, Then $V$ is an irreducible
weight module over $\Vir=\oplus_{i\in\sZ}\C d_i\oplus\C c$.
Consequently, all weight spaces of $V$ are $1$-dimensional.}

{\it Proof}. If $f(t)=p(t)q(t)$ with nonconstant $p(t),q(t)\in\C[t]$
being relatively prime, then we have a natural isomorphism $\varphi:
\Bbb L(f)\rightarrow \Bbb L(p)\oplus \Bbb L(q)$ with
$\varphi(d_0)=(d_0,d_0)$. (Note the different meanings of each
$d_0$.) We consider $\Bbb L(p)\oplus \Bbb L(q)$ in stead of $\Bbb
L(f)$.

\medskip
{\bf Claim 1.} $(d_0,0)$ and $(0,d_0)$ are diagonalizable on $V$.
\medskip

 Suppose that $V=\oplus_{i\in\sZ} V_{i}$,
where $V_{i}=\{v\in V\,|\, d_{0}v=(\a+i)v\}$ and $\a\in\C$. Since
$[(d_0,d_0), (0,d_0)]=0$ and $[(d_0,d_0), (d_0,0)]=0$, there exist a
nonzero $v\in V_i$ such that $(d_0,0)v=\lambda v$ and $(0, d_0)v=\mu
v$. Since $(0,d_0)$ and $(0,d_0)$ are diagonalizable on $L(f)$ and
$V$ is irreducible over $L(f)$, then $(d_0,0)$ and $(0,d_0)$ are
diagonalizable on $V=U(L(f))v$. Claim 1 follows.

Now we suppose eigenvalues of $(d_0,0)$ and $(0,d_0)$ are in
$\lambda+\Z$ and $\mu+\Z$ respectively with $0\le Re(\lambda)<1 \rm
$ {and } $0\le Re(\mu)<1$.

\medskip
{\bf Claim 2.} Either $\Bbb L(p)$ or $\Bbb L(q)$ acts trivially on
$V$.
\medskip

 Suppose on the contrary that both $\Bbb
L(p)$ and $\Bbb L(q)$ act nontrivially on $V$. Considered as an
$\Bbb L(p)$ module, $V$ has a weight space decomposition: $V\cong
\bigoplus_{k\in\sZ}V_{(k, *)}$£¬ where $V_{(k,*)}=\{v\in
V|(d_0,0)v=(\lambda+k)v\}$, and $V_{(k,*)}\neq 0, \forall k\neq
0$. Note that each $V_{(k,*)}$ is an $\Bbb L(q)$ module since
$[\Bbb L(p),\Bbb L(q)]=0$.
\smallskip

If $\Bbb L(q)V_{(k,*)}=0$ for any fixed $k\neq 0$, then $\Bbb L(q)$
acts trivially on $V=U(\Bbb L(p))V_{(k,*)}$. So $\Bbb
L(q)V_{(k,*)}\neq0$. We then have that $V_{(k,*)}\cong
\bigoplus_{j\in\sZ}V_{(k,j)}$, where $V_{(k,j)}=\{v\in
V_{(k,*)}|(0,d_0)v=(\mu+j)v\}$, and $V_{(k,j)}\neq 0, \forall j\neq
0$.
\smallskip

Clearly, all $V_{(k,j)}$ intersects trivially with each other for
different $(k,j)$, and $V_{\lambda+\mu}=\{v\in
V|(d_0,d_0)v=(\lambda+\mu)v\}$ contains
$\bigoplus_{k+j=0}V_{(k,j)}$ which is infinite dimensional.
Contradiction! Thus $\Bbb L(p)$ or $\Bbb L(q)$ acts trivially on
$V$. Claim 2 follows.
\medskip

Now we can write $f(t)=f_1(t)f_2(t)\cdots f_m(t)$ where
$f_k(t)=(t-a_k)^{n_k}$ and $a_k\ne0$ are pairwise distinct complex
numbers. Applying   Claim 2 and using induction on $m$, we can
easily prove that $V$ can be viewed as a uniformly bounded
irreducible weight module over $\Bbb L(f_k)\simeq \Vir\otimes
\Big(\C[t]/(x^{n_k}\C[x])\Big)$ for some $k$. Then by Theorem 4.2,
$V$ is an irreducible weight module over $\Vir=\oplus_{i\in\sZ}\C
d_i\oplus\C c$. Consequently, all weight spaces of $V$ are
$1$-dimensional and $(\Vir\otimes x\C[x])V=0$.\qed

Now we are ready to prove

{\bf Theorem 4.4.}\  {\it If $V$ is a uniformly bounded nontrivial
irreducible module  over $\LL$, then all weight spaces of $V$ are
$1$-dimensional.}

{\it Proof}. From Theorem 4.3 we need only to show that there
exists $R(t)\in\C[t]$ with nonzero constant term such that
$(\Vir\otimes R(t)\C[t^{\pm1}])V=0$.

Suppose that $V=\oplus_{i\in\sZ} V_{i}$, where $V_{i}=\{v\in
V\,|\, d_{0}v=(\a+i)v\}$ and $\a\in\C$. We may assume that
$V_i\ne0$ for any $i\in\Z$ with $\a+i\ne0$.  Take a basis
$\{v^{(1)},v^{(2)},\cdots ,v^{(p)}\}$ of $V_i$ for a fixed
$i\in\Z$. For any $j\neq 0$, we define a linear map:
$$\varphi_j:\,\,\, \C[t^{\pm1}]\,\rightarrow\,
V_{i+j}=V_{i+j}\oplus V_{i+j}\oplus\cdots\oplus V_{i+j},$$
 by $\varphi_j(P(t))=\Big((d_{j}\otimes P(t)) v^{(1)},(d_{j}\otimes
P(t)) v^{(2)},\cdots,(d_{j}\otimes P(t)) v^{(p)}\Big).$ Clearly
$(d_{j}\otimes P(t))V_i=0$ for any $P(t)\in\ker \varphi_j$.
\smallskip

For any $v\in V_i$,  $P(t)\in\ker \varphi_j$  and
$Q(t)\in\C[t^{\pm1}]$, we have that
$$j(d_j\otimes P(t)Q(t))v=[d_0\otimes Q(t), d_j\otimes
P(t)]v$$ $$=(d_0\otimes Q(t))(d_j\otimes P(t))v-(d_j\otimes
P(t))(d_0\otimes Q(t))v=0.$$ Then $\ker \varphi_j$ is an ideal of
$\C[t^{\pm1}]$. Since $\C[t^{\pm1}]$ is a principal ideal domain,
there exists $P_j(t)\in\C[t]$ with nonzero constant term and $\deg
P_j(t)\leq \dim V_i\dim V_{i+j}$ such that $\ker \varphi_j=\langle
P_j(t)\rangle$. Since $V$ is uniformly bounded, we may suppose that
$\deg P_j(t)\leq q$ for some $q\in\N$.
\smallskip

Since $[d_1\otimes P_1(t), d_2\otimes P_2(t)]v=(d_3\otimes
P_1(t)P_2(t))v=0,\,\forall \, v\in V_i$,  then
$P_3(t)|P_1(t)P_2(t)$. Inductively, we can easily get that
$P_j(t)|P_1^{j-2}(t)P_2(t),\,\forall j\geq 3$. Similarly, we also
get that $P_{-j}(t)|P_{-1}^{j-2}(t)P_{-2}(t),\,\forall j\geq 3$.
\smallskip

Now let $R(t)=P_1^{q}(t)P_2(t)P_{-1}^{q}(t)P_{-2}(t).$ Since $\deg
P_j(t)\leq q$ we can easily see that $P_j(t)|R(t)$ for any nonzero
$j\in\Z$. For any $Q(t)\in\langle R(t)\rangle $, we can write
$Q(t)=Q_1(t)Q_{-1}(t)=Q_2(t)Q_{-2}(t)$ where $Q_j(t)\in\langle
P_j(t)\rangle ,\,j\in\{\pm1,\pm2\}$. Then consider $[d_{-1}\otimes
Q_{-1}(t),d_1\otimes Q_1(t)]=2d_0\otimes Q(t)$ and $[d_{-2}\otimes
Q_{-2}(t),d_2\otimes Q_2(t)]=4d_0\otimes Q(t)-(c\otimes Q(t))/2,$
acting on $V_i$, we get that $(d_0\otimes Q(t))V_i=0$ and
$(c\otimes Q(t))V_i=0$. That is, $(\Vir\otimes \langle
R(t)\rangle) V_{i}=0$.
\smallskip

Since $V=U(\LL)V_{i}$, using the PBW basis it is not hard to
deduce that $(\Vir\otimes \langle R(t)\rangle) V=0$. This
completes the proof.\qed

Combining Theorems 3.1 and 4.4 we obtain

{\bf Theorem 4.5.} {\it  Let $V$ be an irreducible Harish-Chandra
$\LL-$module. Then $V$ is either a highest weight module,  a
lowest weight module or a module of the intermediate series.}

From Corollary 3.2 and  Theorem 4.2 we have

 {\bf Theorem 4.6.} {\it Let $V$ be an irreducible Harish-Chandra
module over the truncated Virasoro algebra $\Vir\otimes
\C[t]/\langle t^n\rangle$ for any $n\in\N$. Then $V$ is either a
highest weight module, a lowest weight module or module $V'(a,b)$
over the Virasoro algebra $\Vir={\rm span}\{d_i,c\,|i\in\Z\}$ with
$\Big(\Vir\otimes t\C[t]\Big)V=0$.}

This corollary with $n=2$ classifies all irreducible
Harish-Chandra module over the algebra $W(2,2)$ studied in [ZD].
The classification problem for $W(2,2)$ was also discussed in
[LiZ].

The method in this section (mainly Theorems 4.1 and 4.2) does not
apply to the twisted Heisenberg-Virasoro algebra, that is,
computations in Sect.4 in [LZ3] cannot be simplified.

 \vskip .5cm
\par
\cl{{\bf \S5. Evaluation modules of the intermediate series}}
\par

In this section, we shall determine all indecomposable weight
$\LL$-modules with all weight spaces $1$-dimensional which we have
called modules of the intermediate series.

We shall denote the 1-dimensional trivial module over Vir or $\LL$
by $T$.

{\bf Theorem 5.1.}\ {\it Modules of the intermediate series over
$\LL$ are all evaluation modules of the intermediate series:
$V(a,b)(e), A_a(e), B_a(e), V'(0,0)(e)$ for some $a, b, e\in\C$
with $e\ne0$.}

{\it Proof.} Let $V$ be an $\LL$-module of the intermediate
series. Considered as a $\Vir$-module, $V$ must be one of
$V(\a,\b)$, $A_a$, $B_b$, $V'(0,0)\oplus T$ or $V'(0,0)$. If
$V\simeq A_a$ we set $\a=0$ and $\b=1$;  if $V\simeq   B_b$, or
$V'(0,0)\oplus T$ or $V'(0,0)$,  we set $\a=\b=0$.

We have the weight space decomposition $V=\bigoplus_{i\in\sZ}V_i$,
where $V_i=\{u\in V|d_0u=(\a+i)u\}$. Assume that $\{v_i\in
V_i|i\in\Z\}$ is the standard basis of $V$ defined as before,
except for the case $V\cong V'(0,0)$ where we assign $v_0=0$.

Now suppose that $(d_i\o t^m) v_k=f(i,k,m)v_{i+k}$. Applying
$[d_i\o t^m, d_j\o
t^n]=((j-i)d_{i+j}+\frac{1}{12}\d_{i+j,0}(i^3-i)c)\o t^{m+n}$ to
$v_k$, we obtain that
$$f(i,k+j,m)f(j,k,n)-f(i,k,m)f(j,k+i,n)=(j-i)f(i+j,k,m+n),
 \eqno (5.1)$$
for $i,j,m,n\in\Z$ with $i+j\ne0$ $\Big($and
$(\a+k)(\a+k+i)(\a+k+j)(\a+k+i+j)\ne0$  if $V= V'(0,0)$$\Big)$.

{\bf Case 1.} $V\cong V(\a,\b)$ as $\Vir$-module.

 Since $d_iv_k=(\a+k+i\b)v_{i+k}$, we
see that $f(i,k,0)=\a+k+i\b,$ for all  $i,k\in \Z$.

Taking $j=i$ and $m=0$ in (5.1), we obtain that
$f(i,k+i,0)f(i,k,n)-f(i,k,0)f(i,k+i,n)=0,$ i.e.,
$$ f(i,k,n):(\a+k+i\b)=f(i,k+i,n):(\a+k+i+i\b).
 \eqno (5.2)$$

We denote $k_i=-\a-i\b$.

From (5.2) with $i=1$, we obtain that,   for the case
$\a+\b\in\Z$,
$$\frac{f(1,k_1-i,n)}{\a+\b+k_1-i}=F_1^-(n),\,\,\,{\rm{and}}
 \,\,\,\frac{f(1,k_1+i,n)}{\a+\b+k_1+i}=F_1^+(n),\,\,\,\f i\in\N;\eqno(5.3)$$
$$f(1,k_1,n)=0; \eqno(5.4)$$
where $F_1^+(n)$ and $F_1^-(n)$ are scalars depending on $n$. If
$\a+\b\notin\Z$, we have simpler formula
$$F_1^-(n)=F_1^+(n) = \frac{f(1,k,n)}{\a+\b+k} , \f k\in \Z.\eqno (5.5)$$
It is clear that $F_1^{\pm}(0)=1$.

From (5.2) with $i=-1$, we obtain that,  for the case
$\a-\b\in\Z$,
$$F_{-1}^-(n)= \frac{f(-1,k_{-1}-i,n)}{\a-\b+k_{-1}-i},\,\,\,{\rm{and}}
 \,\,\,
F_{-1}^+(n)=\frac{f(-1,k_{-1}+i,n)}{\a-\b+k_{-1}+i},\,\,\,\f
i\in\N;\eqno (5.6)$$
$$f(-1,k_{-1},n)=0 ;\eqno
(5.7)$$ where $F_{-1}^\pm(n)$   are scalars depending on $n$. If
$\a-\b\notin\Z$, we have simpler formula
$$F_{-1}^-(n)=F_{-1}^+(n) = \frac{f(-1,k,n)}{\a-\b+k}
, \f k\in \Z.\eqno (5.8)$$ Clearly, $F_{-1}^{\pm}(0)=1$.

From (5.2) with $i=2$, we obtain that,   for the case
$\a+2\b\in\Z$,
$$F_2^-(n)=\frac{f(2,k_2-2i,n)}{\a+2\b+k_2-2i} ,\,\,\,{\rm{and}}
 \,\,\,F_2^+(n)=\frac{f(2,k_2+2i,n)}{\a+2\b+k_2+2i},\,\,\,\f
i\in\N; \eqno (5.9)$$
$$f(2,k_2,n)=0;\eqno (5.10)$$
$$G_{2}(n)=\frac{f(2,k_2+2i+1,n)}{\a+2\b+k_2+2i+1},\,\,\,\f
i\in\Z;\eqno (5.11)$$ where $F_2^\pm(n)$ and $G_2(n)$ are scalars
depending on $n$. In the case that $\a+2\b\notin\Z$,  we fix
$k_2=0$ for the sake of consistence although it is an abuse of
$k_2$. Thus we have simpler formulas
$$F_2^-(n)= F_2^+(n)=\frac{f(2,k_2+2i,n)}{\a+2\b+k_2+2i},\,\,\,\f
i\in\Z; \eqno (5.12)$$
$$G_{2}(n)=\frac{f(2,k_2+2i+1,n)}{\a+2\b+k_2+2i+1},\,\,\,\f
i\in\Z;\eqno (5.13)$$
 Note that $F_2^{\pm}(0)=G_2(0)=1$.

From (5.2) with $i=-2$, we obtain that,   for the case
$\a-2\b\in\Z$,
$$F_{-2}^-(n)= \frac{f(-2,k_{-2}-2i,n)}{\a-2\b+k_{-2}-2i},\,\,\,{\rm{and}}
 \,\,\,F_{-2}^+(n)=\frac{f(-2,k_{-2}+2i,n)}{\a-2\b+k_{-2}+2i},\,\,\,\f
i\in\N; \eqno (5.14)$$
$$f(-2,k_{-2},n)=0;\eqno (5.15)$$
$$G_{-2}(n)=\frac{f(-2,k_{-2}+2i+1,n)}{\a-2\b+k_{-2}+2i+1}
 ,\,\,\,\f
i\in\Z;\eqno (5.16)$$ where $F_{-2}^\pm(n)$ and $G_{-2}(n)$ are
scalars depending on $n$. In the case of $\a-2\b\notin\Z$, we fix
$k_2=0$ for the sake of consistence although it is an abuse of
$k_2$. Thus we have
$$F_{-2}^-(n)=F_{-2}^+(n) =\frac{f(-2,k_{-2}-2i,n)}{\a-2\b+k_{-2}-2i},\,\,\,\f
i\in\Z; \eqno (5.17)$$
$$G_{-2}(n)=\frac{f(-2,k_{-2}+2i+1,n)}{\a-2\b+k_{-2}+2i+1}
 ,\,\,\,\f
i\in\Z;\eqno (5.18)$$ Note that $F_{-2}^{\pm}(0)=G_{-2}(0)=1$.

Taking  $i=0$ in (5.1), we see that
$$f(j,k,n)(f(0,k+j,m)-f(0,k,m))=jf(j,k,m+n). \eqno (5.19)$$

For $k\gg0$, using (5.3) and (5.5) to (5.19) with $j=1$ we obtain
that
$$F^+_1(m+n)=F^+_1(n)(f(0,k+1,m)-f(0,k,m)),\f m,n\in\Z\eqno (5.20)$$
Letting $m=-n$, we see that $1=F^+_1(n)(f(0,k+1,-n)-f(0,k,-n))$,
which implies $$F^+_1(n)\ne0,\forall n\in\Z.\eqno (5.21)$$ From
(5.20) we also deduce that
$$F^+_1(n)=(f(0,k+1,1)-f(0,k,1))^n, \f n\in \Z, k\gg0.\eqno(5.22)$$

Similar to (5.20-22), we can deduce that
$$F^\pm_{\pm1}(n)\ne0, F^\pm_{\pm2}(n)\ne0,G _{\pm2}(n)\ne0;\eqno(5.23)$$
$$F^-_1(n)=(f(0,k+1,1)-f(0,k,1))^n, \f n\in \Z, k<\hskip -4pt<0;\eqno(5.24)$$
$$F^+_{-1}(n)=(f(0,k,1)-f(0,k-1,1))^n, \f n\in \Z, k\gg0;\eqno(5.25)$$
$$F^-_{-1}(n)=(f(0,k,1)-f(0,k-1,1))^n, \f n\in \Z, k<\hskip -4pt<0;\eqno(5.26)$$
$$F^+_2(n)=(\frac{1}{2}(f(0,k+2,1)-f(0,k,1)))^n, \f n\in \Z, k\gg0,k\equiv k_2(\mod 2);\eqno(5.27)$$
$$F^-_2(n)=(\frac{1}{2}(f(0,k+2,1)-f(0,k,1)))^n, \f n\in \Z, k<\hskip -4pt<0,k\equiv k_2(\mod  2);\eqno(5.28)$$
$$F^+_{-2}(n)=(\frac{1}{2}(f(0,k ,1)-f(0,k-2,1)))^n, \f n\in \Z, k\gg0;k\equiv k_{-2}(\mod 2);\eqno(5.29)$$
$$F^-_{-2}(n)=(\frac{1}{2}(f(0,k,1)-f(0,k-2,1)))^n, \f n\in \Z, k<\hskip -4pt<0,k\equiv k_{-2}(\mod  2);\eqno(5.30)$$
$$G_2(n)=(\frac{1}{2}(f(0,k+2,1)-f(0,k,1)))^n, \f n\in \Z,k\equiv k_2+1(\mod 2);\eqno(5.31)$$
$$G_{-2}(n)=(\frac{1}{2}(f(0,k,1)-f(0,k-2,1)))^n, \f n\in \Z,k\equiv k_2+1(\mod 2).\eqno(5.32)$$
From (5.24) and (5.25) we see that
$f(0,k,1)-f(0,k-1,1)=f(0,k+j,1)-f(0,k-1+j,1)$ and
$f(0,-k-j,1)-f(0,-k-j-1,1)=f(0,-k,1)-f(0,-k-1,1)$ for $k\gg0$ and
$j\in\N$. Using these formulas and (5.22)-(5.32) we deduce that
$$F_1^+(n)=F_{-1}^+(n)=F_2^+(n)=F_{-2}^+(n)=G_2(n)$$ $$=G_{-2}(n)
=F_1^-(n)=F_{-1}^-(n)=F_2^-(n)=F_{-2}^-(n)=e^n, \f n\in
Z\eqno(5.33)$$ for some nonzero $e\in\C$.

Now we have obtained that $$(d_i\o t^m)v_k=(\a+i\b+k)e^mv_{i+k},\f
k,m\in\Z, i=\pm1,\pm2.$$

Since the Lie algebra $\LL$ is generated by $d_i\o t^m$ for all $
m\in\Z, i=\pm1,\pm2$,  $V$ has to be the evaluation module
$V(\a,\b)(e)$.

{\bf Case 2.} $V\cong B_b$ as a $\Vir$-module.

In this case, we have $\a=\b=0$, and $f(i,k,0)=(\a+i\b+k)$ except
that $f(i,-i,0)=-i(i+b)$. Thus the formula (5.1) still holds. We
replace (5.2)  with
$$ f(i,k,n):f(i,k,0)=f(i,k+i,n):f(i,k+i,0),
 \eqno (5.2')$$

Similar to (5.3)-(5.18), we have
$$F_1^-(n)= \frac{f(1,-i-1,n)}{-i-1}=\frac{f(1,-1,n)}{-(1+b)},\f i\in\N;\eqno (5.34)$$
$$f(1,0,n)=0;\eqno (5.35)$$
$$F_1^+(n)=\frac{f(1,i,n)}{i},\f i\in\N;\eqno (5.36)
$$
$$F_{-1}^-(n)=\frac{f(-1,-i,n)}{-i},\f i\in\N;\eqno (5.37)$$
$$f(-1,0,n)=0;\eqno (5.38)$$
$$F_{-1}^+(n)=\frac{f(-1,1,n)}{b-1}=\frac{f(-1,i+1,n)}{i+1},\f i\in\N;\eqno (5.39)
$$
$$F_2^-(n)=\frac{f(2,-2i-2,n)}{-2i-2}=\frac{f(2,-2,n)}{-2(b+2)},\f i\in\N;\eqno (5.40)$$
$$f(2,0,n)=0;\eqno (5.41)$$
$$F_2^+(n)=\frac{f(2,2i,n)}{2i},\f i\in\N;\eqno (5.42)$$
$$G_{2}(n)=\frac{f(2,2i+1,n)}{2i+1} ,\f i\in\Z;\eqno (5.43) $$
$$F_{-2}^-(n)=\frac{f(-2,-2i,n)}{-2i},\f i\in\N;\eqno (5.44)$$
$$f(-2,0,n)=0;\eqno (5.45)$$
$$F_{-2}^+(n)=\frac{f(-2,2,n)}{2(b-2)}=\frac{f(-2,2i+2,n)}{2i+2},\f i\in\N;\eqno (5.46)$$
$$G_{-2}(n)=\frac{f(-2,2i+1,n)}{2i+1},\f i\in\Z.\eqno (5.47)$$
Note that if any denominator involving $b$ in the above  formulas
is $0$,  then we shall next deduce that the corresponding
numerator must be $0$.

With the same argument after (5.18) in Case 1, we can deduce that
$V\cong B_b(e)$ for some nonzero $e\in\C$.

{\bf Case 3.} $V\cong V'(0,0)$ as a $\Vir$-module.

In this case, the coefficients $f(i,k,n)$ are defined only for
$k(k+i)\ne0$, and we have   $f(i,k,0)=k$. The formula (5.1) holds
for $k(k+i)(k+j)(k+i+j)\ne0$, and (5.2) holds for
$k(k+i)(k+2i)\ne0$.

In this case, similar to (5.34)-(5.47), we have
$$F_1^-(n)= \frac{f(1,-i-1,n)}{-i-1},\f i\in\N;$$
$$F_1^+(n)=\frac{f(1,i,n)}{i},\f i\in\N;
$$
$$F_{-1}^-(n)=\frac{f(-1,-i,n)}{-i},\f i\in\N;$$
$$F_{-1}^+(n)=\frac{f(-1,i+1,n)}{i+1},\f i\in\N;
$$
$$F_2^-(n)=\frac{f(2,-2i-2,n)}{-2i-2},\f i\in\N;$$
$$F_2^+(n)=\frac{f(2,2i,n)}{2i},\f i\in\N;$$
$$G_{2}(n)=\frac{f(2,2i+1,n)}{2i+1} ,\f i\in\Z; $$
$$F_{-2}^-(n)=\frac{f(-2,-2i,n)}{-2i},\f i\in\N;$$
$$F_{-2}^+(n)=\frac{f(-2,2i+2,n)}{2i+2},\f i\in\N;$$
$$G_{-2}(n)=\frac{f(-2,2i+1,n)}{2i+1},\f i\in\Z.$$

With the same argument after (5.18) in Case 1, we can deduce that
$V\cong V'(0,0)(e)$  for some nonzero $e\in\C$.

{\bf Case 4.} $V\cong V'(0,0)\oplus T$ as a $\Vir$-module.

This case is almost the same as Case 2 with some modifications.
The modifications are the following. First note that
$f(i,-i,0)=f(i,0,0)=0$ for all $i\in\Z$. Then in (5.34), (5.39),
(5.40), (5.46), delete the fractions involving $b$, and add the
identities that all the corresponding numerators equal $0$. Follow
the proof in Case 2,
 we deduce that $V\cong (V'(0,0)\oplus T)(e)$
for some nonzero $e\in\C$. Thus $V\cong V'(0,0)(e)\oplus T$.

{\bf Case 5.} $V\cong A_a$ as a $\Vir$-module.

Consider the graded dual module $V^*$ of $V$.  As a $\Vir$-module,
 $V^*=\oplus_{i\in\sZ}\C v^*_i$ is
isomorphic to $B_a$. From Case 2, we know that $V^*\simeq B_a(e)$
for some nonzero $e\in\C$. Hence $V\cong A_a(e)$.

Now we have proved that $V$ is one of evaluation modules of the
intermediate series. \qed

As a direct consequence, using the isomorphism $\Vir\otimes
\Big(\C[t]/t^n\C[t]\Big)\simeq \Vir\otimes
\Big(\C[t^{\pm1}]/(t+1)^n\C[t^{\pm1}]\Big)$ we obtain

 {\bf Theorem 5.2.}\  {\it If $V$ is an irreducible weight module over the truncated Virasoro algebra
 $\Vir\otimes\Big(\C[t]/\langle t^n\rangle\Big)$ with all weight spaces 1-dimensional, where $n\in\N$, then
$\Big[\Vir\otimes t\Big(\C[t]/\langle t^n\rangle\Big)\Big]\cdot
V=0$. Consequently, $V$ is an irreducible  module over
$\Vir=\oplus_{i\in\sZ}\C d_i\oplus\C c$.}

We summarize the established results into the following

{\bf Theorem 5.3.}\  {\it (a) If $V$ is an irreducible weight
module over $\LL$,  then
 $V$ is a highest weight module, a lowest weight module, or
 $V'(a,b)(e)$ for some $a,b,e\in\C$ with $c\ne0$.

(b) If $V$ is an irreducible weight module over the truncated
Virasoro algebra
 $\Vir\otimes\Big(\C[t]/\langle t^n\rangle\Big)$, where $n\in\N$, then
 $V$ is a highest weight module, a lowest weight module or  an irreducible  module over
$\Vir=\oplus_{i\in\sZ}\C d_i\oplus\C c$ with all weight spaces
$1$-dimensional and $\Big[\Vir\otimes t\Big(\C[t]/\langle
t^n\rangle\Big)\Big]\cdot V=0$.}

\vskip .5cm

\par
\cl{{\bf \S6. Highest weight modules over $\LL$}}
\par
\vskip .2cm

 Let us study the highest weight modules $V(\varphi)$ define in (2.2).
 Generally, not all weight spaces $V_{\lambda-i}$
 of $V(\varphi)$ are  finite-dimensional. The following theorem
 answers this question.

\medskip
 {\bf Theorem 6.1.} {\it The weight module $V(\varphi)$ is a Harish-Chandra module if
and only if there is a polynomial $P(t)\in \C[t]$ with nonzero
constant term such that $$\varphi(d_0\otimes
t^kP(t))=0,\,\,\,\,\varphi(c\otimes t^kP(t))=0\f k\in \Z,$$ in which
case, $(d_i\o \la P(t)\ra )V(\varphi)=0, \f  i\in\Z$, where $\la
P(t)\ra $ is the ideal of $\C[t^{\pm}]$ generated by $P(t)$.
Consequently, $V(\varphi)$ can be viewed as a module over $\Vir\o
(\C[t^{\pm}]/\la P(t)\ra )$ if it is a Harish-Chandra module. }

{\it Proof.} First suppose that $V(\varphi)$ is a Harish-Chandra
module. Then there is some $P(t)\in \C[t]$ with nonzero constant
term such that $(d_{-2}\o t^sP(t))v_0=0$ for some $s\in\Z$. Apply
$d_2\o t^k$ and $(d_1\o t^k)(d_1\o 1)$ to it for any $k\in\Z$, we
deduce that:
$$
0=(d_2\o t^k)(d_{-2}\o t^sP(t))v_0=((-4d_0+\frac{1}{2}c)\o
t^{k+s}P(t))v,
$$
$$
0=(d_1\o t^k)(d_1\o 1)(d_{-2}\o t^sP(t))v_0=(6d_0\o
t^{k+s}P(t))v_0.
$$
Hence, $(d_0\o t^{k}P(t))v_0=0$ and $(c\o t^{k}P(t))v_0=0$, i.e.,
$\varphi(d_0\otimes t^kP(t))=0$ and $\varphi(c\otimes t^kP(t))=0$
for all $k\in\Z$.

Conversely, suppose that there is a polynomial $P(t)\in \C[t]$ with
nonzero constant term such that $\varphi(d_0\otimes t^kP(t))=0$ and
$\varphi(c\otimes t^kP(t))=0$, i.e., $(d_0\otimes t^kP(t))v_0=0$ and
$(c\otimes t^kP(t))v_0=0$ for all $k\in \Z$.

We shall prove that $(d_{n}\otimes t^kP(t))v_0=0$ for any $n,k\in
\Z$ by induction on $n$. This is clear for $n\geq0$. Now suppose
it true for $n\ge-(i-1)$ where $ i\geq 1$. Consider
$(d_{-i}\otimes t^kP(t))v_0, \forall k\in\Z$.  For any $s\in \Z$,
we deduce that
$$
(d_1\o t^s)(d_{-i}\otimes t^kP(t))v_0=((-i-1)d_{-i+1}\o
t^{k+s}P(t))v_0=0,
$$
$$
(d_2\o t^s)(d_{-i}\otimes
t^kP(t))v_0=(((-i-2)d_{-i+2}+\frac{1}{2}\d_{i,2}c)\o
t^{k+s}P(t))v_0=0.
$$
If $(d_{-i}\otimes t^kP(t))v_0\ne0$, it would be a highest weight
vector with lower weight than that of $v_0$, contrary to the
irreducibility of $V(\varphi)$. Thus $(d_{n}\otimes t^kP(t))v_0=0$
for any $n,k\in \Z$.

Using PBW  theorem and by induction on $i\ge0$ we can easily deduce
that $(d_{n}\otimes t^kP(t))V(\varphi)_{\lambda-i}=0$ for all $n,
k\in\Z$. Hence $(d_{n}\otimes t^kP(t))V(\varphi) =0$ for all $n,
k\in\Z$.

Now $V(\varphi)$ can be viewed as an irreducible module over $\Vir\o
\Big(\C[t^{\pm}]/\la P(t)\ra\Big )$. Since the weight spaces of
$\Vir\o \Big(\C[t^{\pm}]/\la P(t)\ra \Big)$ are all finite
dimensional, so do those of $V(\varphi)$ by PBW theorem.\qed

\par
It is easy to see that the evaluation module $M(\dot c,h;a)$ is a
module $V(\varphi)$ in Theorem 6.1 with $P(t)=t-a$,
$\varphi(d_0)=h$ and $\varphi(c)=\dot c$.  Thus we shall call
irreducible Harish-Chandra modules $V(\varphi)$ in Theorem 6.1
with $P(t)=(t-a)^n$ {\bf generalized evaluation highest weight
modules} over $\LL$.

For convenience, for any $P(t)\in\C[t^{\pm1}]$ let us define $$
\hom_P(\LL_0,\C)=\{\varphi\in \hom(\LL_0,\C)|\varphi(d_0\otimes
t^kP(t))=\varphi(c\otimes t^kP(t))=0 \f k\in\Z\}.$$ Since
$\hom_P(\LL_0,\C)\simeq \hom (\Vir_0\otimes \C[t^{\pm1}]/\la
P(t)\ra,\C)$, we shall identify them. If $P=P_1P_2..P_r$ with
$(P_i,P_j)=1$ for all $i\ne j$. Then we have the algebra isomorphism
$ \C[t^{\pm}]/\la P(t)\ra \cong\bigoplus_{i=1}^{r}(\C[t^{\pm}]/\la
P_i\ra ). $ Thus, we have the corresponding Lie algebra isomorphism
$$\Vir\o\C[t^{\pm}]/\la P(t)\ra
\simeq\bigoplus_{i=1}^{r}(\Vir\o(\C[t^{\pm}]/\la P_i\ra
)).\eqno(6.1)$$ We shall identify these two algebras by the
natural map, that is, we can write any element $x\otimes f(t)$,
where $x\in\Vir$ and $f(t)\in\C[t^{\pm}]$, as
$$
x\otimes f(t)=(x\otimes f(t),x\otimes f(t),\cdots,x\otimes f(t)).
$$
Actually the precise expression should be
$$
x\otimes f(t)(\mod (P(t))=(x\otimes f(t)(\mod P_1),x\otimes
f(t)(\mod P_2),\cdots,x\otimes f(t)(\mod P_r)).
$$
 Note that, in this notation, we have
 $$
[(x_1,x_2,\cdots,x_r),(y_1,y_2,\cdots,y_r)]=([x_1,y_1],
\cdots,[x_r,y_r]),
 $$
where $x_{1},x_2,...,x_r,y_1...y_r\in \LL$. For simplicity, we
shall write $\Vir\o\Big(\C[t^{\pm}]/\la P(t)\ra\Big)$ and
$\Vir_0\o\Big(\C[t^{\pm}]/\la P(t)\ra\Big)$ as $\LL/\la P(t)\ra$
and $\LL_0/\la P(t)\ra$ respectively.

Then we have the following natural identifications:
$$\LL/\la P(t)\ra=\oplus _{i=1}^r\left( \LL/\la
P_i\ra\right),\eqno(6.2)$$
$$\hom_P(\LL_0,\C)=\oplus _{i=1}^r
\hom_{P_i}(\LL_0,\C).\eqno(6.3)$$

Thus for any $\varphi\in \hom_P(\LL_0,\C)$ there exist
$\varphi_i\in \hom_{P_i}(\LL_0,\C)$ such that $\varphi=\sum
_{i=1}^r\varphi_i$. Then we have the tensor product  $\LL$-module
$V(\varphi_1)\otimes V(\varphi_2)\otimes\cdots \otimes
V(\varphi_r)$. If we use the identification (6.2), then
$x=(x_1,\cdots, x_r)$ with $x_i\in \LL/\la P_i\ra$, and
$$x(v_1\otimes v_2\otimes\cdots \otimes v_r)= x_1v_1\otimes
v_2\otimes\cdots \otimes v_r\hskip 3cm
$$
$$+v_1\otimes x_2v_2\otimes\cdots \otimes v_r+\cdots+ v_1\otimes
v_2\otimes\cdots \otimes x_rv_r.$$

{\bf Proposition 6.2.} Suppose $P=P_1P_2\in\C[t^{\pm1}]$ with
$(P_1, P_2)=1$, and $\psi=\psi_1+\psi_2\in \hom_{P}(\LL_0,\C)$
with $\psi_i\in \hom_{P_i}(\LL_0,\C)$. Then $V(\psi)\simeq
V(\psi_1)\otimes V(\psi_2)$.

{\it Proof.} Since $V(\psi_1)$ and $V(\psi_2)$ are weight modules,
so is $V(\psi_1)\o V(\psi_2)$. We denote the highest weight vector
of $V(\psi_i)$ by $v^{(i)}$. We need only prove that any nonzero
homogeneous element in $V(\psi_1)\o V(\psi_2)$ can generate
$v^{(1)}\o v^{(2)}$.

For any fixed  nonzero homogeneous element $u\in V(\psi_1)\o
V(\psi_2)$, let $$u=\sum_{j=1}^n u_j^{(1)}\o u_j^{(2)},
u_j^{(i)}\in V^{(i)}_{\l^{(i)}_j}$$ with
$\l^{(1)}_j+\l^{(2)}_j=\l^{(1)}_{j'}+\l^{(2)}_{j'}$. We may
suppose that $n$ is minimal in such expressions. Thus we may also
assume that $\l^{(1)}_1= \cdots =\l^{(1)}_s<\l^{(1)}_{s+1}\leq
\cdots \leq \l^{(1)}_n$. Then $u^{(1)}_1, u^{(1)}_2, \cdots ,
u^{(1)}_s$ are linearly independent and $u^{(2)}_1, u^{(2)}_2,
\cdots , u^{(2)}_s$ are linearly independent.

Because of the identification (6.2), instead of using $\LL$ (or
$\LL/\la P\ra$) we shall use $\LL/\la P_1\ra\oplus \LL/\la
P_2\ra$. Since $V(\psi_1)$ is an irreducible weight module over
$\LL/\la P_1\ra$, then there is a homogeneous element $x\in U(\LL
/\la P_1\ra)$, the universal enveloping algebra, such that
$xu^{(1)}_1=v^{(1)}$. Thus   $x(u^{(1)}_1\o
u^{(2)}_1)=(xu^{(1)}_1)\o u^{(2)}_1=v^{(1)}\o u^{(2)}_1$ since
$x(u^{(2)}_j)=0$. It is clear that $x(u^{(1)}_j\o u^{(2)}_j)\in
(\C v^{(1)})\o u^{(2)}_j$ for $ 1\leq j\leq s$ and that
$x(u^{(1)}_j\o u^{(2)}_j)=0$ for $ s+1\leq j\leq n$. Noting that
$u^{(2)}_1, u^{(2)}_2, \cdots , u^{(2)}_s$ are linearly
independent, we know that $0\neq xu=v^{(1)}\otimes w\in (\C
v^{(1)})\o V^{(2)}$. By similar discussions, we can find $x'\in
U(\LL /\la P_2\ra)$ such that $x'w= v^{(2)}$, i.e.,
$x'xu=v^{(1)}\o x'w=v^{(1)}\o v^{(2)}$. This completes the proof.
\qed

Directly from Proposition 6.2 we have

 {\bf Corollary 6.3.} {\it Let $\varphi\in \hom_{P}(\LL_0,\C)$
where $P(t)=\prod _{i=1}^{r}(t-a_i)^{n_i}$ with $n_i\in\N, a_i\neq
a_j\in\C\setminus\{0\}, \forall i\neq j$. Denote
$P_i=(t-a_i)^{n_i}$. Then there exist $\varphi_i\in
\hom_{P_i}(\LL_0,\C)$ such that $\varphi=\varphi_1+ \varphi_2+
\cdots +\varphi_r$ and $V(\varphi)\simeq V(\varphi_1)\otimes
V(\varphi_2)\otimes\cdots \otimes V(\varphi_r)$.}

Now applying Corollary 6.3  we can deduce the following:

{\bf Theorem 6.4.} {\it Any irreducible highest weight module with
finite dimensional weight spaces over $\Vir\o \C[t^{\pm1}]$ is a
tensor product of some generalized evaluation modules.}

Remark that the irreducibility of the Verma modules over the
truncated  algebras \break $\LL/\langle (t-a)^n\rangle$ for any
$n\in\N$ and $a\in\C$ is determined in [Wi].

In the rest of this section, we shall study the necessary and
sufficient conditions for  the Verma module $\bar V(\varphi)$ over
$\LL$ to be irreducible. The answer   is the following:

\par
{\bf Theorem 6.5.}  {\it The Verma module $\bar V(\varphi)$ is not
irreducible if and only if there  exists a polynomial $P(t)\in
\C[t]$ with nonzero constant term such that
$$\varphi(d_0\otimes t^kP(t))=0 \f
k\in \Z.$$}
\par
We see from Theorems 6.1 and 6.5 that we cannot always get
Harish-Chandra modules $V(\varphi)$ from not irreducible Verma
modules $V(\bar\varphi)$.

Before proving  Theorem 6.5, we need some preparations on the
universal enveloping algebra $U=U(\LL_-)$ of $\LL_-$.

Let $\SS$ be the set of finite sequences of integers
$(i_1,i_2,\cdots,i_r)$. We first define a partial ordering $\succ$
on the set $\SS$: $(i_1,i_2,\cdots,i_r) \succ
(j_1,j_2,\cdots,j_s)$ if and only if $i_1=j_1$, $i_2,=j_2, \cdots,
i_{k-1}=j_{k-1}$ and $i_k>j_k$.

We have the obvious meaning for $\succeq, \preceq$ and $\prec$.

We fix a PBW basis $\B$ for $U(\LL_-)$ consisting of the following
elements:
$$(d_{-i_1}\o t^{j_1})(d_{-i_2}\o t^{j_2})...(d_{-i_r}\o t^{j_r});\,\,\,
 (i_s,j_s)\in \N\times\Z, r\in\Z_+$$
 where $(i_s,j_s)\succeq (i_{s+1},j_{s+1})$ for all
$s=1,2,...,r-1$. We call $r$ the {\bf height} of the element
$(d_{-i_1}\o t^{j_1})(d_{-i_2}\o t^{j_2})...(d_{-i_r}\o t^{j_r})$,
which is denoted by $\ht((d_{-i_1}\o t^{j_1})(d_{-i_2}\o
t^{j_2})...(d_{-i_r}\o t^{j_r}))$. We now define an ordering on
$\B$ as follows:
$$(d_{-i_1}\o t^{j_1})(d_{-i_2}\o t^{j_2})...(d_{-i_r}\o t^{j_r})\succ
(d_{-i'_1}\o t^{j'_1})(d_{-i'_2}\o t^{j'_2})...(d_{-i'_s}\o
t^{j'_s})$$ if $(r,i_1,i_2,\cdots,i_r,j_1,j_2,\cdots,j_r)\succ
(s,i'_1,i'_2,\cdots,i'_s,j'_1,j'_2,\cdots,j'_s)$. Then $\succ$ is
a total ordering on $\B$.

Any nonzero $X\in U(\LL_-)$  can be uniquely written  as a linear
combination of elements in $\B$: $X=\Sigma_{i=1}^m a_iX_i$ where $
0\neq a_i\in\C, X_i\in {\B}$ and $X_1\succ X_2\succ\cdots\succ
X_m.$ We define the {\bf height of} $X$ as $\ht(X)=\ht(X_1)$, and
the {\bf highest term} of $X$ as $\hm(X)=a_1X_1$. For convenience,
we define $\ht(0)=-1$ and $\hm(0)=0$.

It is clear that $\B v_0:=\{Xv_0\,\,|\,\,X\in\B\}$ is a basis for
the Verma module $\bar V(\varphi)$ where $v_0$ is again the
highest weight vector of $\bar V(\varphi)$. We define
$$ \ht(Xv_0):=\ht(X),\hskip .5cm \hm(Xv_0):=\hm(X)v_0, \hskip 1cm \forall\hskip 2pt
X \in U(\LL_-).$$

 We need to define the last notation
$$U^r_s=\{X\in U(\LL_-)\,\,|\,\,[d_0,X]=sX, \ht(X)\le r\}.$$ It is
easy to see that $U^r_{-s}U^{r'}_{-s'}\subset U^{r+r'}_{-s-s'}$
for all $r,r',s,s'\in\N$.

Also we have $(d_1\o t^k)U^r_{-s} v_0 \subset  U^r_{-s+1} v_0$ for
all $r,s\in\Z_+$ and $k\in\Z$ (where we have regarded $U^r_{k}=0$
for $k>0$) and $(d_{-i_1}\o t^{j_1})(d_{-i_2}\o
t^{j_2})...(d_{-i_r}\o t^{j_r})\in U_{-i_1-i_2-\cdots-i_r}^r$ for
any $i_1,i_2,...,i_r\in\N$ and $j_1,j_2,...,j_r\in\Z$.

\par
Now we are ready to give
\par
 {\it Proof of Theorem 6.5.} Suppose there   exists a polynomial
$P(t)\in \C[t]$ with nonzero constant term such that
$\varphi(d_0\otimes t^kP(t))=0$ for all $ k\in \Z.$ We can easily
deduce that, in $V(\varphi)$,  $(d_{-1}\otimes t^kP(t))v_0=0$ for
all $k\in \Z.$ Thus $\bar V(\varphi)$ is not irreducible.

Now suppose that there  does not exist a polynomial $P(t)\in \C[t]$
with nonzero constant term such that $\varphi(d_0\otimes t^kP(t))=0$
for all $k\in \Z.$ We want to prove that $\bar V(\varphi)$ is
irreducible. In fact, we need only to show that $\bar
 V(\varphi)=V(\varphi)$, i.e., $\bar
 V(\varphi)_{-n}=V(\varphi)_{-n}$ for all $n\in \N$.
 We shall do this by induction on $n$. The statement for $n=0$ is trivial.

Now consider $n=1$. If $ \bar
 V(\varphi)_{-1}\neq  V(\varphi)_{-1}$, then there is some $P(t)\in \C[t^{\pm}]$ such
 that $(d_{-1}\o P(t))v_0=0$ in $V(\varphi)$. Then,  in $V(\varphi)$,
 $0=(d_1\o t^k)(d_{-1}\o P(t))v_0=
 [d_1\o t^k,d_{-1}\o P(t)]v_0=-2(d_0\o t^kP(t))v_0=-2\varphi(d_0\o t^kP(t))v_0 $, i.e.,
 $\varphi(d_0\o t^kP(t)) =0$ for all $k\in\Z$, a contradiction. Thus $ \bar
 V(\varphi)_{-1}=V(\varphi)_{-1}$.

Now we consider the case $n>1$ and suppose that $ \bar V_{-k}=
V_{-k}$ for all $0\le k<n$.

 It suffices to prove that $Xv_{0}\neq 0$ in $V(\varphi)$ for
 any  nonzero $X\in U(\LL_-)_{-n}$. We write
 $X=\Sigma_{i=1}^m a_iX_i$ where $ 0\neq a_i\in\C, X_i\in
{{\B}}$ with $X_1\succ X_2\succ\cdots\succ X_m.$ We are going to
show that $Xv_0\neq0$. We break up the proof into two different
cases.

{\bf{Case 1.}}   $\ht(X)<n$.

Suppose $X_1=(d_{-i_1}\o t^{j_1})(d_{-i_2}\o
t^{j_2})...(d_{-i_r}\o t^{j_r}) (d_{-1}\o t^{j_{r+1}})(d_{-1}\o
t^{j_{r+2}})...(d_{-1}\o t^{j_{r+s}}) \in\B$ with $r>0$, and
$i_r\geq 2$. Since $(d_1\o t^k)(Xv_0)\in V_{-(n-1)}=\bar
V_{-(n-1)}$,  by the inductive hypothesis we know that, for
sufficiently large $k\in \N$,
$$\hm((d_1\o t^k)(Xv_0))=\hm([d_1\o t^k,a_1X_{1}]v_0)$$
$$=-pa_1(i_r+1)(d_{-i_1}\o t^{j_1})(d_{-i_2}\o t^{j_2})...(d_{-i_r+1}\o t^{j_r+k})
(d_{-1}\o t^{j_{r+1}})(d_{-1}\o t^{j_{r+2}})...(d_{-1}\o
t^{j_{r+s}})v_0$$ $$\neq 0,$$ where $p$ is the number of $q$ such
that $(i_q,j_q)=(i_r,j_r)$. Then $(d_1\o t^k)(Xv_0)\neq0$ in
$V(\varphi)$, yielding $Xv_0\neq0$ in $V(\varphi)$.

{\bf{Case 2.}} $\ht(X)=n$.

To the contrary in this case, we assume that $Xv_0=0$ in
$V(\varphi)$. There exist $r,s$ such that $\ht(X_i)=n, 1\leq i\leq
r$, $\ht(X_i)=n-1, r+1\leq i\leq s$ and $\ht(X_i)\leq n-2, s+1\leq
i\leq m$.

For $1\leq i\leq r$, each $X_i$ is of the form $X_i=(d_{-1}\o
t^{j^{(i)}_{1}})(d_{-1}\o t^{j^{(i)}_{2}})...(d_{-1}\o
t^{j^{(i)}_{n}})$. For any $k\in\Z$  we compute
$$(d_1\o t^k)(X_iv_0)=[d_1\o t^k,(d_{-1}\o t^{j^{(i)}_{1}})(d_{-1}\o t^{j^{(i)}_{2}})
...(d_{-1}\o t^{j^{(i)}_{n}})]v_0$$
$$
=-2\sum_{p=1}^{n} (d_{-1}\o t^{j^{(i)}_1})(d_{-1}\o
t^{j^{(i)}_2})... (d_{-1}\o t^{j^{(i)}_{p-1}})(d_0\o
t^{j^{(i)}_p+k})(d_{-1}\o t^{j^{(i)}_{p+1}})...(d_{-1}\o
t^{j^{(i)}_n})v_0
$$
$$
=-2\sum_{p=1}^{n} (d_{-1}\o t^{j^{(i)}_1})(d_{-1}\o
t^{j^{(i)}_2})... \widehat{(d_{-1}\o t^{j^{(i)}_{p}})}...
(d_{-1}\o t^{j^{(i)}_n})(d_0\o t^{j^{(i)}_p+k})v_0
$$
$$+2\sum_{p=1}^{n}\sum_{q=p+1}^n
(d_{-1}\o t^{j^{(i)}_1})... \widehat{(d_{-1}\o
t^{j^{(i)}_{p}})}... (d_{-1}\o t^{j^{(i)}_{q-1}})(d_{-1}\o
t^{j^{(i)}_p+j^{(i)}_q+k}) (d_{-1}\o
t^{j^{(i)}_{q+1}})...(d_{-1}\o t^{j^{(i)}_n})v_0$$
$$\equiv-2\sum_{p=1}^{n}
(d_{-1}\o t^{j^{(i)}_1})(d_{-1}\o t^{j^{(i)}_2})...
\widehat{(d_{-1}\o t^{j^{(i)}_{p}})}... (d_{-1}\o
t^{j^{(i)}_n})\varphi(d_0\o t^{j^{(i)}_p+k})v_0$$
$$+2\sum_{p=1}^{n}\sum_{q=p+1}^n
(d_{-1}\o t^{j^{(i)}_p+j^{(i)}_q+k})(d_{-1}\o
t^{j^{(i)}_1})...\widehat{(d_{-1}\o t^{j^{(i)}_p})}...
\widehat{(d_{-1}\o t^{j^{(i)}_{q}})}...(d_{-1}\o t^{j^{(i)}_n})v_0$$
$$\mod(U^{(n-2)}_{-n+1}v_0), \eqno(6.4)$$ where the $\hat{ }$ means
the factor is omitted.

For $r+1\leq i\leq s$, each $X_i$ is of the form $X_i=(d_{-2}\o
t^{j^{(i)}_{1}})(d_{-1}\o t^{j^{(i)}_{2}})...(d_{-1}\o
t^{j^{(i)}_{n-1}})$. For any $k\in\Z$  we compute
$$(d_1\o t^k)(X_iv_0)=
[d_1\o t^k,(d_{-2}\o t^{j^{(i)}_{1}})(d_{-1}\o
t^{j^{(i)}_{2}})...(d_{-1}\o t^{j^{(i)}_{n-1}})]v_0$$
$$\equiv
-3(d_{-1}\o t^{j^{(i)}_1+k})(d_{-1}\o t^{j^{(i)}_{2}})...
(d_{-1}\o t^{j^{(i)}_{n-1}})v_0
\,\,\,\mod(U^{(n-2)}_{-n+1}v_0).\eqno(6.5)$$
 Using (6.4),
(6.5) and the fact that $(d_1\o t^k)(X_iv_0)\in U^{(n-2)}_{-n+1}$
for all $s+1\leq i\leq m$, we obtain that, for all $k\in\Z$,
$$0=(d_1\o t^k)(Xv_0)=[d_1\o t^k, X]v_0\hskip 7cm$$
$$\equiv-2\sum_{i=1}^{r}a_i\sum_{p=1}^{n}
(d_{-1}\o t^{j^{(i)}_1})(d_{-1}\o t^{j^{(i)}_2})...
\widehat{(d_{-1}\o t^{j^{(i)}_{p}})}... (d_{-1}\o
t^{j^{(i)}_n})\varphi(d_0\o t^{j^{(i)}_p+k})v_0$$
$$+2\sum_{i=1}^{r}a_i\sum_{p=1}^{n}\sum_{q=p+1}^n
(d_{-1}\o t^{j^{(i)}_p+j^{(i)}_q+k})(d_{-1}\o
t^{j^{(i)}_1})...\widehat{(d_{-1}\o t^{j^{(i)}_p})}...
\widehat{(d_{-1}\o t^{j^{(i)}_{q}})}...(d_{-1}\o
t^{j^{(i)}_n})v_0$$
$$
-3\sum_{i=r+1}^{s}a_i(d_{-1}\o t^{j^{(i)}_1+k})(d_{-1}\o
t^{j^{(i)}_{2}})... (d_{-1}\o t^{j^{(i)}_{n-1}})v_0
\,\,\,\mod(U^{(n-2)}_{-n+1}v_0)$$
$$\equiv-2\sum_{i=1}^{r}a_i\sum_{p=1}^{n}
(d_{-1}\o t^{j^{(i)}_1})(d_{-1}\o t^{j^{(i)}_2})...
\widehat{(d_{-1}\o t^{j^{(i)}_{p}})}... (d_{-1}\o
t^{j^{(i)}_n})\varphi(d_0\o t^{j^{(i)}_p+k})v_0$$
$$+\sum_{{\underline{l}}}
g_{{\underline{l}}}(k)(d_{-1}\o t^{l_1+k})(d_{-1}\o
t^{l_2})...(d_{-1}\o
t^{l_{n-1}})v_0\,\,\,\mod(U^{(n-2)}_{-n+1}v_0), \eqno(6.6)$$ where
the last summation is a finite one over
${\underline{l}}=(l_1,l_2,...,l_{n-1})$ with $l_2\geq l_3...\geq
l_{n-1}$, and $g_{{\underline{l}}}(k)$ are polynomials in $k$.
Denote $$x_1=-2\sum_{i=1}^{r}a_i\sum_{p=1}^{n} (d_{-1}\o
t^{j^{(i)}_1})(d_{-1}\o t^{j^{(i)}_2})... \widehat{(d_{-1}\o
t^{j^{(i)}_{p}})}... (d_{-1}\o t^{j^{(i)}_n})\varphi(d_0\o
t^{j^{(i)}_p+k})v_0$$
$$x_2=\sum_{{\underline{l}}}
g_{{\underline{l}}}(k)(d_{-1}\o t^{l_1+k})(d_{-1}\o
t^{l_2})...(d_{-1}\o t^{l_{n-1}})v_0.$$

Now for any sufficiently large $k\in \Z$, let $$R=\{ (d_{-1}\o
t^{j^{(i)}_1})(d_{-1}\o t^{j^{(i)}_2})... \widehat{(d_{-1}\o
t^{j^{(i)}_{p}})}... (d_{-1}\o t^{j^{(i)}_n})v_0\,\,|\,\,
i=1,2,...,r \}$$ which is the set of all possible basis elements
appearing in the expression of $x_1\in V_{-(n-1)}=\bar V_{-(n-1)}
$, and let
$$T=\{(d_{-1}\o t^{l_1+k})(d_{-1}\o t^{l_2})...(d_{-1}\o t^{l_{n-1}})v_0\}$$
which is the set of all  possible basis elements appearing in the
expression of  $x_2\in V_{-(n-1)}=\bar V_{-(n-1)} $. By inductive
hypothesis we know that $R$ and $T$ are linearly independent in
the vector space $(U^{(n-1)}_{-n+1}v_0)/(U^{(n-2)}_{-n+1}v_0)$ for
sufficiently large $k$. Thus $x_1=0=x_2$ for sufficiently large
$k$. Consequently $g_{{\underline{l}}}(k)=0$ for all $k$, i.e.,
$x_1=x_2=0$ for all $k\in\Z$.

In terms of linear combination of $\B v_0$, the coefficient of
$(d_{-1}\o t^{j^{(1)}_1})(d_{-1}\o t^{j^{(1)}_2})... (d_{-1}\o
t^{j^{(1)}_{n-1}})v_0$ in the expression of $x_1$ for all $k\in
\Z$ should be $0$. On the other hand, it is
$$-2\sum_{i\in I}p_ia_i
\varphi(d_0\o t^{j^{(i)}_n+k})=\varphi(d_0\o t^k(-2\sum_{i\in I}(
p_ia_i t^{j^{(i)}_n})))$$ where $I=\{1\leq i\leq
r|(j^{(i)}_1,j^{(i)}_2,...,j^{(i)}_{n-1})=(j^{(1)}_1,j^{(1)}_2,...,j^{(1)}_{n-1})\}$,
$p_1$ is the number of $q$ such that $j^{(1)}_n=j^{(1)}_q$ and other
$p_i=1$. If we denote $P(t_2)=-2\sum_{i\in I}( p_ia_i
t^{j^{(i)}_n})=\sum_{j=1}^{h}h_j t^j$,  then $P(t_2)\ne0$. Hence
$(d_{0}\o t^kP(t))v_0=0, \f k\in\Z$, i.e., $\varphi(d_0\otimes
t^kP(t))=0, \forall k\in \Z$, a contradiction. Therefore $Xv_0\neq0$
in this case.

This completes the proof of the Theorem 6.5.\qed

Let us recall  exp-polynomial functions defined in [BZ]. The
algebra of  {\bf exp-polynomial functions} in the integer variable
$i$ is the algebra generated as an algebra by functions
 $f(i):\Z\rightarrow \C$ with $f(i)=i$ and
$f(i)=a^{i}$ for various constants $a\in
\C^{\star}=\C\backslash\{0\}.$ Use a well-known combinatoric
formula (see [RZ]), we can restate Theorems 6.1 and 6.4 as
follows.

\medskip
 {\bf Theorem 6.6.} {\it (a) The weight module $V(\varphi)$ is a Harish-Chandra module if
and only if there are exp-polynomial functions  $f(k), g(k)$ such
that $$\varphi(d_0\otimes t^k )=f(k),\,\,\,\,\varphi(c\otimes t^k
)=g(k)\f k\in \Z.$$

(b) The Verma module $\bar V(\varphi)$ is reducible if and only if
there is an exp-polynomial function  $f(k)$ such that
$\varphi(d_0\otimes t^k )=f(k)\f k\in \Z.$}

We conclude this section by giving an example of generalized
evaluation module over $\LL$.

{\bf Example.} Let us take $\varphi(d_0(k))=(-1)^kk$ and
$\varphi(c(k))=0$ for all $k\in\Z$. Then $V(\varphi)$ is a module
over $\Vir\otimes \Big(\C[t^{\pm1}]/(t+1)^2\C[t^{\pm1}]\Big)$.
Clearly, $V(\varphi)$ is not an evaluation module. It is a
generalized evaluation module with $P(t)=(t+1)^2$.  Although
$\varphi(d_0)=\varphi(c)=0$, you can easily see that $d_{-1}(0)v_0$
and $d_{-1}(1)v_0$ are linearly independent in $V(\varphi)$.
Actually, from Corollary 7.11 in [Wi], this module is an irreducible
Verma module over the truncated algebra $\Vir\otimes
\Big(\C[t^{\pm1}]/(t+1)^2\C[t^{\pm1}]\Big)\simeq \Vir\otimes
\Big(\C[t]/t^2\C[t]\Big)$.

 \vskip 20pt
 \cl{{\bf References}}
 \par
\begin{description}

 \item{[A]} Robert B. Ash, {\it Basic abstract algebra}, Dover Publications,
 Mineola Inc., NY, 2007.

 \item{[BGLZ]}  P. Batra, X. Guo, R. Lu  and  K. Zhao, Highest weight modules
over the pre-exp-polynomial algebras, preprint, February, 2008.

 \item{[BZ]}  Y. Billig and K. Zhao,  Weight
modules over exp-polynomial Lie algebras, {\it J. Pure Appl.
Algebra}, Vol.191, 23-42(2004).

\item{[DL]} C. Dong and J. Lepowsky, {\it Generalized vertex
algebras and relative vertex operators}, Progress in Math. Vol.
112, Birkhauser, Boston, 1993.

\item{[FMS]} P. Di Francesco, P. Mathieu and D. Senechal,  {\it
Conformal Field Theory}, Springer, New York, 1997.

\item{[H]} J.~E.~Humphreys, {\it{Introduction to Lie algebras and
representation theory}}, 2nd ed. (rev.), Springer-Verlag, Berlin,
1972.

 \item{[IKU]} T. Inami, H. Kanno and T. Ueno,
Higher-dimensional WZW model on Kahler manifold and toroidal Lie
algebra, {\it Mod. Phys. Lett.} A 12, 2757-2764 (1997).

\item{[IKUX]} T. Inami, H. Kanno, T. Ueno and C.-S. Xiong,
Two-toroidal Lie algebra as current algebra of four-dimensional
Kahler WZW model, {\it Phys.Lett.} B 399, 97-104 (1997).

 \item{[K]} V.G. Kac, {\it Infinite-dimensional Lie algebras},
Third edition, Cambridge University Press, Cambridge, 1990.

\item{[K1]} I. Kaplansky, The Virasoro algebra, {\it Comm. Math.
Phys.}, 86(1982), no.1, 49-54.

\item{[KR]} V.Kac and A. Raina, {\it Bombay lectures on highest
weight representations of infinite dimensional Lie algebras},
World Sci., Singapore, 1987.

\item{[LiZ]} D. Liu and L. Zhu, Classification of Harish-Chandra
modules over the $W$ algebra $W(2,2)$, ARXIV: 0801.2601v2.

\item{[LZ]} R. Lu and K. Zhao, Verma modules over quantum torus
 Lie algebras, to appear in {\it Canadian J. Math.}

\item{[LZ2]} R. Lu and  K. Zhao, Classification of irreducible
weight modules over higher rank Virasoro algebras, {\it Advances
in Math.}, Vol.201(2), 630-656(2006).

\item{[LZ3]} R. Lu and K. Zhao, Classification of irreducible
weight modules over the twisted Heisenberg-Virasoro algebras,
ARXIV: 0510190.

 \item{[M]}  I.G. Macdonald,  Affine
root systems and Dedekind's $\eta $-function,
 {\it Invent. Math.}15(1972), 91--143.

\item{[Ma]} O. Mathieu, Classification of Harish-Chandra modules
over the Virasoro algebra, {\it Invent. Math.} 107(1992), 225-234.

\item{[MP]} C. Martin and A. Piard,  Nonbounded indecomposable
admissible modules over the Virasoro algebra, {\it Lett. Math.
Phys.} 23(1991), no. 4, 319--324.

\item{[MoP]} R. V. Moody and  A. Pianzola,  {\it Lie algebras with
triangular decompositions}, Canad. Math. Soc., Ser. Mono. Adv.
Texts, A Wiley-Interscience Publication, John Wiley \& Sons Inc.,
New York, 1995.

 \item{[MZ]} V. Mazorchuk and  K. Zhao,   Classification of
simple weight Virasoro modules with a finite-dimensional weight
space, {\it J. Algebra}, Vol.307, 209--214(2007).

\item{[RZ]} S. Eswara Rao and K. Zhao, Highest weight irreducible
representations of quantum tori, {\it Mathematical Research
Letters}, Vol.11, No.5-6, 615-628(2004).

\item{[SZ]} Y. Su and K. Zhao, Generalized Virasoro and
super-Virasoro algebras and modules of intermediate series, {\it
J. Algebra}, 252(2002), no.1, 1-19.

\item{[W]} X. Wang,  Evaluation modules of Lie algebras with a
triangular decomposition,  {\it Adv. Math.}(China), 33(2004), no.6,
685-690.

\item{[Wi]} B. J. Wilson, Highest weight theory for truncated
current Lie algebras, ARXIV: 0705.1203.

\item{[ZD]}  W. Zhang and C. Dong, W-algebra $W(2,2)$ and the
vertex operator algebra $L(1/2,0)\otimes L(1/2,0)$, ARXIV: 071104.
\end{description}

Department of Mathematics, Zhengzhou university, Zhengzhou 450001,
Henan, P. R. China. Email: guoxq@amss.ac.cm

Department of Mathematics, Suzhou university, Suzhou 215006,
Jiangsu, P. R. China. Email: rencail@amss.ac.cn

Department of Mathematics, Wilfrid Laurier University, Waterloo,
ON, Canada N2L 3C5, and Academy of Mathematics and System
Sciences, Chinese Academy of Sciences, Beijing 100190, P. R.
China. Email: kzhao@wlu.ca

\end{document}